\documentclass[]{article}
\usepackage{graphicx}
\usepackage{caption}
\usepackage{algorithmic}
\usepackage{algorithm}
\usepackage{amsmath}
\usepackage{amssymb}
\usepackage{amsfonts}
\usepackage{amsthm}
\usepackage{fullpage}
\usepackage{url}
\usepackage{color}
\usepackage{authblk} 
\usepackage{float}
\usepackage{mathabx} 
\usepackage{subfig}

\def\N{{\mathbb{N}}}
\def\R{{\mathbb{R}}}

\def\dt{{\textrm{d}t}}

\def\IMF{{\textrm{IMF}}}
\def\IMFs{{\textrm{IMFs}}}

\begin{document}

\title{A novel algorithm for the decomposition of non-stationary multidimensional and multivariate signals}
\author{Roberto Cavassi\thanks{DISIM, Universit\'a degli Studi dell'Aquila, L'Aquila, Italy ({\tt roberto.cavassi@graduate.univaq.it})}\and Antonio Cicone\thanks{DISIM, Universit\'a degli Studi dell'Aquila, L'Aquila, Italy\\Istituto di Astrofisica e Planetologia Spaziali, INAF, Via del Fosso del Cavaliere 100, 00133,Rome, Italy\\Istituto Nazionale di Geofisica e Vulcanologia, Via di Vigna Murata 605, 00143, Rome, Italy ({\tt antonio.cicone@univaq.it})}
\and Enza Pellegrino\thanks{DIIIE, Universit\'a degli Studi dell'Aquila, L'Aquila, Italy({\tt enza.pellegrino@univaq.it})}
\and Haomin Zhou\thanks{School of Mathematics, Georgia Institute of Technology, Atlanta, GA, U.S.A. ({\tt hmzhou@math.gatech.edu})}}

\maketitle

\begin{abstract}
The decomposition of a signal is a fundamental tool in many fields of research, including signal processing, geophysics, astrophysics, engineering, medicine, and many more. By breaking down complex signals into simpler oscillatory components we can enhance the understanding and processing of the data, unveiling hidden information contained in them. Traditional methods, such as Fourier analysis and wavelet transforms, which are effective in handling mono-dimensional stationary signals struggle with non-stationary data sets and they require, this is the case of the wavelet, the selection of predefined basis functions. In contrast, the Empirical Mode Decomposition (EMD) method and its variants, such as Iterative Filtering (IF), have emerged as effective nonlinear approaches, adapting to signals without any need for a priori assumptions. To accelerate these methods, the Fast Iterative Filtering (FIF) algorithm was developed, and further extensions, such as Multivariate FIF (MvFIF) and Multidimensional FIF (FIF2), have been proposed to handle higher-dimensional data. 

In this work, we introduce the Multidimensional and Multivariate Fast Iterative Filtering (MdMvFIF) technique, an innovative method that extends FIF to handle data that vary simultaneously in space and time. This new algorithm is capable of extracting Intrinsic Mode Functions (IMFs) from complex signals that vary in both space and time, overcoming limitations found in prior methods. The potentiality of the proposed method is demonstrated through applications to artificial and real-life signals, highlighting its versatility and effectiveness in decomposing multidimensional and multivariate nonstationary signals. The MdMvFIF method offers a powerful tool for advanced signal analysis across many scientific and engineering disciplines.
\end{abstract}

%

\section{Introduction}\label{sec:Intro}

Signal decomposition plays an important role in applied mathematics and signal processing, and it has many applications in data analysis, image processing, machine learning, geophysics, communications, and many other fields. Decomposing complicated signals into simpler oscillatory components can help in understanding the underlying structures and enables efficient data representation, noise reduction, and feature extraction. Customary decomposition techniques like Fourier analysis and wavelet transforms effectively process one-dimensional signals and methods that can handle multidimensional and multivariate signals are increasingly important as real-world data grows in complication. However, these methods have limitations. Fourier-based techniques assume the stationarity of the signal, whereas wavelet-based methods require the a priori selection of a proper basis to be used. Furthermore, both approaches are linear in nature. 

To address these limitations, in 1998 the pioneering Empirical Mode Decomposition method was published \cite{huang1998empirical}. EMD separates signals into their oscillatory components, known as Intrinsic Mode Functions (IMFs), by iteratively subtracting local averages defined by the envelopes connecting maxima and minima of the signal. While effective in many applications, EMD and derived algorithms suffer from several limitations, including their reliance on envelope interpolation and the lack of a rigorous mathematical framework. For these reasons in the last two decades many alternative nonlinear methods have been developed. Among them, the only based on iteration, and therefore not requiring any a priori assumption on the signal under investigation, like the number of IMFs contained in the signal, or the basis to be used in the decomposition, is the Iterative Filtering (IF) method.

The Iterative Filtering (IF) algorithm's structure is based on EMD's ones, with the key difference in the way the signal's moving average is computed. IF computes the moving average $\mathcal{M}(s)(x)$ as the convolution of the signal $s$ with a compactly supported filter or window $w$, instead of relying on the average of two envelopes, like in the EMD
\begin{equation}\label{eq:Mov_Average}
	\mathcal{M}(s)(x)=\int_{-L}^{L} s(x+t)w(t)\dt
\end{equation}
Here, $L$ represents the filter length, or half-support of the filter function $w$. Unlike traditional approaches, IF determines $L$ adaptively based on the signal itself. For instance, as suggested in \cite{lin2009iterative,cicone2016adaptive}, $L$ can be computed using the relative distances between subsequent extrema of the signal. This self-adaptive approach makes IF inherently nonlinear, enabling it to adjust to the signal's characteristics.
Building on IF, the Fast Iterative Filtering (FIF) algorithm was developed to accelerate the decomposition process using the Fast Fourier Transform (FFT). FIF preserves the adaptive nature of IF while significantly reducing computational costs, making it suitable for high-dimensional data. In addition to its speed, FIF introduces the ability to select filter properties, such as the Fokker-Planck filter, which ensures smoothness and compact support.

In recent years, FIF has been extended to handle increasingly complex datasets, leading to the development of the Multivariate Fast Iterative Filtering (MvFIF) \cite{cicone2022multivariate}, the Multidimensional Iterative Filtering (MIF) \cite{cicone2017multidimensional} and Multidimensional Fast Iterative Filtering (FIF2) \cite{sfarra2022maximizing} algorithms. These extensions enable the decomposition of multichannel signals and multidimensional data, respectively, adapting the principles of FIF to higher dimensional and temporal domains. 

In this work, we introduce the Multidimensional and Multivariate Fast Iterative Filtering (MdMvFIF) technique which extends FIF to handle datasets that vary, at the same time, across space and time, and we present its application to handle artificial and real-life signals. To the best of our knowledge this algorithm is the first of its kind. There have been few algorithms proposed so far in the literature for the decomposition of nonstationary signals defined in higher dimensions, possibly including also time, like the Fast and adaptive empirical mode decomposition for multidimensional, multivariate signals (FA-MVEMD) \cite{thirumalaisamy2018fast} and the Fast empirical mode decomposition method for multi-dimensional signals based on serialization (Serial-EMD) \cite{zhang2021serial}. However, these methods handle time like if it were an extra spatial dimension and not as a dimension per se. This is fine as far as the time and space variations are separated. However, in real life signals this is usually not the case and space and time variations are tangled. Furthermore, if we want to remove different range of frequencies from time and space it will be impossible to do it with these algorithms. This is a common problem in applications. Think, for instance, to what shown in \cite{Reagan2025}. In that work the data sets need to be filtered in high frequency in time and low frequency in space or vice-versa. This task cannot be accomplished with the algorithms proposed so far in the literature.

The rest of this work is organized as follows. We first briefly recall the multidimensional and the multivariate extensions of the Fast Iterative Filtering. Then, in Section 2, we present the proposed approach to handle data that are, at the same time, multidimensional and multivariate. In Section 3 we present a few examples of application to both artificial and real life signals. We derive the conclusions in Section 4.

\subsection{Multidimensional Fast Iterative Filtering}\label{sec:FIF2}

The Multidimensional Fast Iterative Filtering (FIF2) method allows to decompose a $n$ dimensional signal $s\in\R^n$ into simple oscillatory components, called Intrinsic Mode Functions (IMFs), by approximating the higher dimensional moving average of $s$ and iteratively subtracting it. Unlike the Empirical Mode Decomposition (EMD), which relies on maxima and minima envelopes, FIF computes the moving average using a convolution of $s$ with a multidimensional filter $w$, defined as a nonnegative, even, and continuous function compactly supported in $\Omega\subset\R^n$, with
$\int_{\Omega} w(z) d^n z = 1$.

The algorithm assumes the periodicity of $s$ at its boundaries, a limitation that can be addressed by pre-extending the signal to enforce periodicity \cite{stallone2020new}. The filter support size $\Omega$ is calculated based on the signal’s characteristics, such as its extrema count or Fourier spectrum. This signal-dependent determination of $\Omega$ makes FIF2 nonlinear and ensures adaptive decomposition.

The FIF2 method involves two nested loops:
\begin{itemize}
	\item \textbf{Inner Loop:} Extracts an IMF by iteratively applying the filtering process until a stopping criterion is met.
	\item \textbf{Outer Loop:} Applies the same process to the residual signal, producing successive IMFs until the residual has at most one local extremum, classifying it as a trend signal.
\end{itemize}

The computational complexity of FIF2 is $\mathcal{O}(m \log m)$, where $m$ is the maximum among all dimensions of the signal under investigation. The method’s robustness and efficiency are enhanced by choosing a Generalized Fokker-Planck filter \cite{cicone2017multidimensional}. These filters are $C^\infty(\mathbb{R}^n)$, compactly supported, and widely used in applications.

\subsection{Multivariate Fast Iterative Filtering}\label{sec:MvFIF}

The Multivariate Fast Iterative Filtering (MvFIF) method extends FIF to handle $n$-dimensional signals \cite{cicone2022multivariate}. For a signal $s \in \mathbb{R}^n \times \mathbb{R}$, the algorithm computes a single filter length $L$ and applies the 1D FIF to each channel. 

To determine $L$, the signal is treated as a sequence of column vectors $\mathbf{v}(t)$ rotating in $\mathbb{R}^n$. The filter length is derived from the angle of rotation $\theta(t)$ between consecutive vectors:
\begin{equation}\label{eq:theta}
	\theta(t)=\arccos\left(\frac{\textbf{v}(t)}{\left\|\textbf{v}(t)\right\|}\cdot\frac{ \textbf{v}(t-1)}{\left\|\textbf{v}(t-1)\right\|}\right)
\end{equation}

The filter length is set as twice the average distance between extrema in $\theta(t)$. This approach ensures that the decomposition adapts to the highest frequency rotations embedded in the signal. The computational complexity of MvFIF is $\mathcal{O}(nm \log m)$, where $n$ is the number of channels and $m$ is the signal length.

\section{Proposed Approach}\label{sec:MdMvFIF}

Given an $n+1$ dimensional signal $f(\mathbf{v},t)$, that varies over space as $\mathbf{v}\in\R^n$ and time as $t\in I\subset\R$, we want to decompose it into IMFs. We propose in this work what we call the Multidimensional and Multivariate Fast Iterative Filtering (MdMvFIF) algorithm. The idea is to extract an IMF in space and one in time iteratively until the signal does not contain any oscillations anymore both in space and time. In particular, for $f(\mathbf{v},\tilde{t})\in\R^n$, $\tilde{t}\in I$ fixed, we start estimating the filter support $\Omega_t\in\R^n$ to be applied to $f(\mathbf{v},t)$, as done in the FIF2 algorithm. We estimate $\Omega_t$ for every time instant $t$ and then compute 
\begin{equation}\label{eq:Omega}
  \widehat{\Omega}=\min\left\{\Omega_t\right\}
\end{equation}
which will be the support of the filter we apply to all space tensors $f(\mathbf{v},t)$ for every fixed $t\in I$. Following what was done in FIF2, we select a Generalized Fokker-Planck filter $w_s$ \cite{cicone2017multidimensional} as filter function and we impose its support to be $\widehat{\Omega}$. Then, we convolve it with the signal to produce its moving average $\mathcal{M}_s(f)=f \ast w_s$. We then subtract it from the signal and iterate the procedure producing $f_{k+1}=f_k-\mathcal{M}_s (f_k)$, for $k\in\N$, until a stopping criterion is satisfied. 

Regarding the stopping criterion, there are many possibile options. One way of doing it, as mentioned in \cite{lin2009iterative,cicone2014adaptive} for the monodimensional case, is to consider the relative change in $f_{k}$ and discontinue the inner loop at a certain $k=N_s\in\N$ as soon as a prefixed threshold is reached.

The previous steps will produce the first estimate of the first IMF in space, $\IMF_s^{(1)}=f_{N_s}$. To speed up calculations, we can compute convolution in the frequency domain where it becomes a simple product. This is made possible, with low computational complexity, thanks to the FFT method. To give the measure of the acceleration induced by this alternative approach, we can consider that if the signal under investigation is $f\in\R^2\times\R$ and after discretization is a cube of dimension $n\times m\times p$, then the computational complexity of the convolution passes from $\mathcal{O}(p\times (n\times m)^2)$ in time domain to $\mathcal{O}(p \times (n\times m) \log(n\times m))$ if done in frequency domain via FFT. 

After computing the first IMF in space, we can subtract it from the signal to obtain a new signal $f=f-\IMF_s^{(1)}$. We can then compute the signal moving average over time. To do so we apply a multivariate approach, as done in the Multivariate Fast Iterative Filtering (MvFIF) method \cite{cicone2022multivariate}. We consider each signal $f(\mathbf{v},\tilde{t})$ obtained by fixing $\tilde{t}\in I$ as a vector in a $\R^{n+1}$ space, and we measure its rotations $\tilde{\theta}(t)$ as the time passes
\begin{equation}\label{eq:theta_new}
	\tilde{\theta}(t)=\arccos\left(\frac{f(\mathbf{v},\tilde{t}+1)}{\left\|f(\mathbf{v},\tilde{t}+1)\right\|}\cdot\frac{ f(\mathbf{v},\tilde{t})}{\left\|f(\mathbf{v},\tilde{t})\right\|}\right)
\end{equation}
We point out that in this formula, differently from what it was done in the past, the time variability is now taken into account \cite{cicone2022multivariate}.
We can now define filter length $L$ of a Fokker-Planck filter $w_t$ we want to convolve with the signal to obtain its moving average. $L$, which is the half support size of the filter, can be computed as the double average distance between subsequent extrema in the signal $\tilde{\theta}(t)$. We point out that this approach is very natural if we consider a multivariate IMF as a vector in $\R^n$ rotating around the time axis. Computing the double average distance between subsequent extrema in $\tilde{\theta}(t)$ allows us to estimate the average scale of the highest frequency rotations embedded in this given signal.

We can now convolve the computed filter $w_t$ with the signal to produce its moving average $\mathcal{M}_t(f)=f \ast w_t$. We then subtract it from the signal and iterate the procedure producing $f_{k+1}=f_k-\mathcal{M}_t (f_k)$, for $k\in\N$, until a stopping criterion is satisfied. Also in this case we can consider the relative change in $f_{k}$ and discontinue the inner loop at a certain $k=N_t\in\N$ as soon as a prefixed threshold is reached. This is how we produce the first estimate of the first IMF in time, $\IMF_t^{(1)}=f_{N_t}$. If we subtract this other first IMF from the signal we obtain a new signal $f=f-\IMF_t^{(1)}$ and we can repeat the entire previous procedure to extract the subsequent IMFs in space and time. We iterate this procedure until there are no more oscillations left both in space and time.

The pseudo-code of MdMvFIF is given in Algorithm \ref{algo:MdMvFIF} and a Matlab implementation is available online\footnote{\url{http://www.cicone.com}}.

\begin{algorithm}
\caption{\textbf{Multidimensional and Multivariate Fast Iterative Filtering} IMFs = MdMvFIF$(f)$}\label{algo:MdMvFIF}
\begin{algorithmic}
\STATE IMFs = $\left\{\right\}$
\STATE i=0
\WHILE{there are oscillations left in space or time}	  
      \STATE i=i+1
      \STATE compute the spatial filter support $\widehat{\Omega}$ of the space filter function $w_s$
      \STATE compute the 2D filter $w_s$
      \STATE compute DFT of signal $f$ and of the filter $w_s$
      \STATE set $k=0$
      \WHILE{the stopping criterion is not satisfied}
                 \STATE  $ \IMF_s^{(i)} = \textrm{iDFT}\left(I-\textrm{diag}\left(\textrm{DFT}(w_s)\right)\right)^{k}\textrm{DFT}(\textbf{f})$
            \STATE  $k = k+1$
      \ENDWHILE
      \STATE $\IMFs_s = \IMFs_s \,\cup\,  \left\{ \IMF_s^{(i)}\right\}$
      \STATE $f=f-\IMF_s^{(i)}$
      \STATE compute the time filter length $L$ based on $\tilde{\theta}(t)$, computed as in \eqref{eq:theta_new}
      \STATE compute the filter $w_t$
      \STATE compute DFT of signal $f$ and of the filter $w_t$
      \STATE set $k=0$
      \WHILE{the stopping criterion is not satisfied}
                 \STATE  $ \IMF_t^{(i)} = \textrm{iDFT}\left(I-\textrm{diag}\left(\textrm{DFT}(w_t)\right)\right)^{k}\textrm{DFT}(\textbf{f})$
            \STATE  $k = k+1$
      \ENDWHILE
      \STATE $\IMFs_t = \IMFs_t \,\cup\,  \left\{ \IMF_t^{(i)}\right\}$
      \STATE $f=f-\IMF_t^{(i)}$
\ENDWHILE
\STATE $\IMFs_s = \IMFs_s \,\cup\,  \{ f\}$
\end{algorithmic}
\end{algorithm}

When the space and time variations are separated, i.e. in the signal $f(\mathbf{v},t)$ the variation over time $t$ is the same for every point $\mathbf{v}\in\R^n$, we could further accelerate calculations via the convolution of the data tensor with a $n+1$ dimensional filter. We call this algorithm the Space Time Fast Iterative Filtering (ST\_FIF), whose pseudo-code is given in Algorithm \ref{algo:ST_FIF}

\begin{algorithm}
\caption{\textbf{Space Time Fast Iterative Filtering} IMFs = ST\_FIF$(s)$}\label{algo:ST_FIF}
\begin{algorithmic}
\STATE IMFs = $\left\{\right\}$
\WHILE{the stopping criterion is not satisfied}
	  \STATE compute the time filter length $L$ for the filter based on $\tilde{\theta}(t)$, computed as in \eqref{eq:theta}
      \STATE compute the spatial filter support $\Omega$
      \STATE compute the $n+1$ filter $w$ whose support is $\widehat{\Omega}\times [-L,L]$
      \STATE compute DFT of signal $s$ and of the filter $w$
      \STATE set $k=0$
      \WHILE{the stopping criterion is not satisfied}
                 \STATE  $ \textrm{IMF}_i = \textrm{iDFT}\left(I-\textrm{diag}\left(\textrm{DFT}(w)\right)\right)^{k}\textrm{DFT}(\textbf{s})$
            \STATE  $k = k+1$
      \ENDWHILE
      \STATE IMFs = IMFs$\,\cup\,  \left\{ \textrm{IMF}_i\right\}$
      \STATE $s=s-\textrm{IMF}_i$
\ENDWHILE
\STATE IMFs = IMFs$\,\cup\,  \{ s\}$
\end{algorithmic}
\end{algorithm}
In this pseudocode, DFT and iDFT stand for Discrete Fourier Transform and inverse DFT, respectively, whereas $\textrm{IMF}_i$ represents the i-th IMF. 

ST\_FIF produces a decomposition that is somehow equivalent to what is also done in the Fast and adaptive empirical mode decomposition for multidimensional, multivariate signals (FA-MVEMD) \cite{thirumalaisamy2018fast} and the Fast empirical mode decomposition method for multi-dimensional signals based on serialization (Serial-EMD) \cite{zhang2021serial}. Furthermore, these methods have a second limitation, as explained also in the introduction, the IMFs produced contain time and space frequency combined together. If we want to filter different frequency ranges for time and space it will not be possible with this algorithm. For this reason, in this work we focus only on the first algorithm and we show in the next section a few applications to synthetic and real-life signals.

\section{Numerical Examples}\label{sec:Examples}

To show the abilities of the MdMvFIF algorithm in handling data we apply it to two artificial signals and a geophysical data set. For the sake of simplicity, especially in plotting results, in the following we focus on signals that are 2D in space and that vary over time, producing 3D data cubes. However, it is important to underline here that the proposed approach can be applied to handle signals of any dimension in space plus time. 

\subsection{First artificial example}

In this first example, we have a signal containing two non-stationary frequencies active in space and two stationary frequencies active in time. 
\begin{figure}[ht]
    \centering
\includegraphics[width=0.9\linewidth]{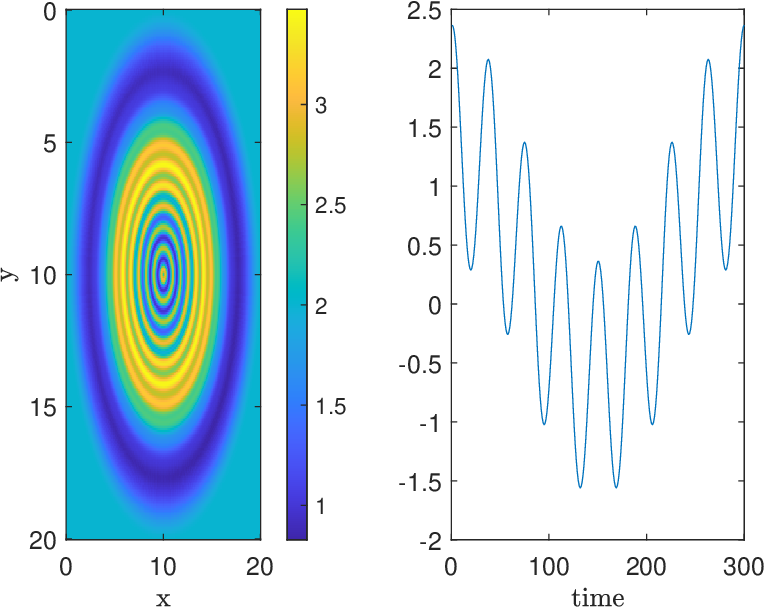}
    \caption{Example 1. Left panel: Signal at time $t=0$. Right panel: Time evolution for position $\mathbf{v}=(150,\,150)$.}
    \label{fig:Ex1_sig}
\end{figure}

\begin{figure}[ht]
    \centering
\includegraphics[width=0.9\linewidth]{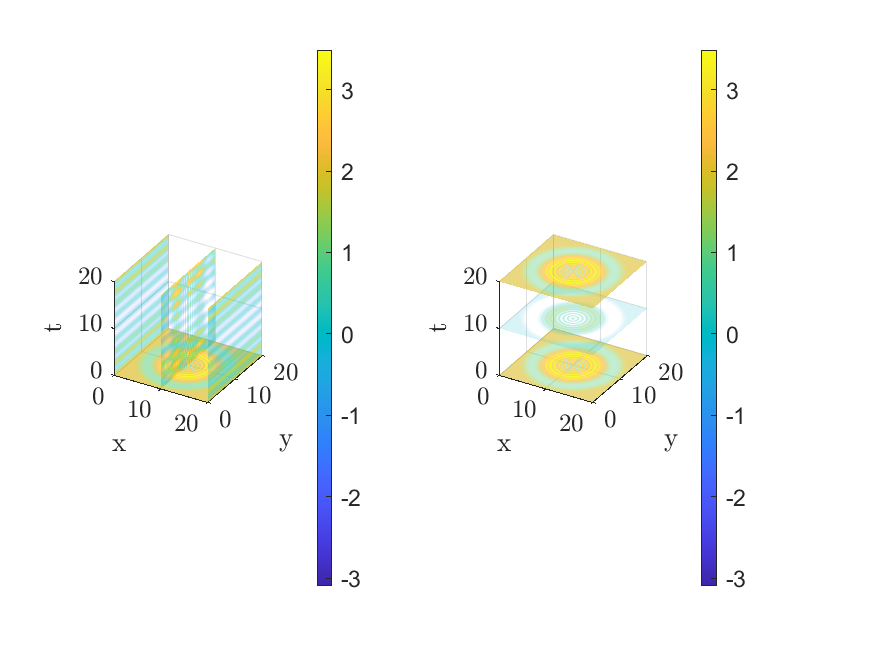}
    \caption{Example 1. Left panel: 3D signal evolution in the second space variable and time. Right panel: 3D signal evolution in space at different time stamps.}
    \label{fig:Ex1_sig2}
\end{figure}

The data are shown in Figures \ref{fig:Ex1_sig} and \ref{fig:Ex1_sig2}. 

\begin{figure}[ht]
    \centering
\includegraphics[width=0.45\linewidth]{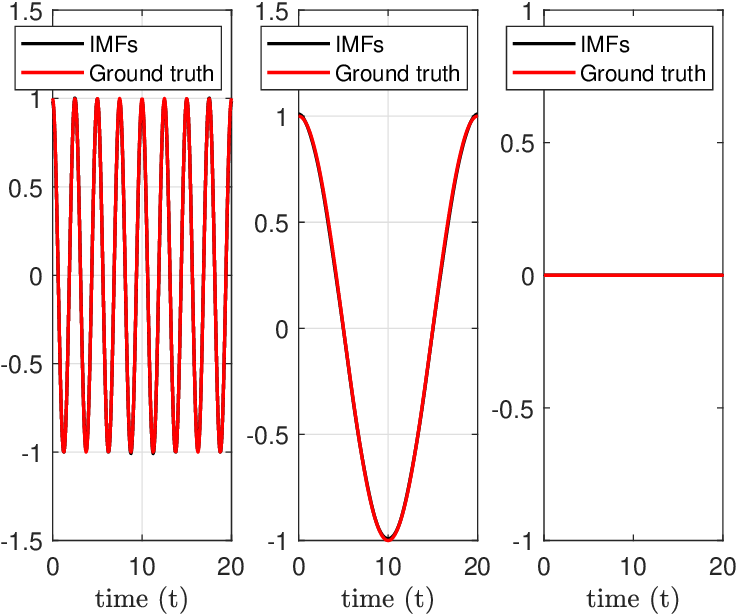}
\includegraphics[width=0.45\linewidth]{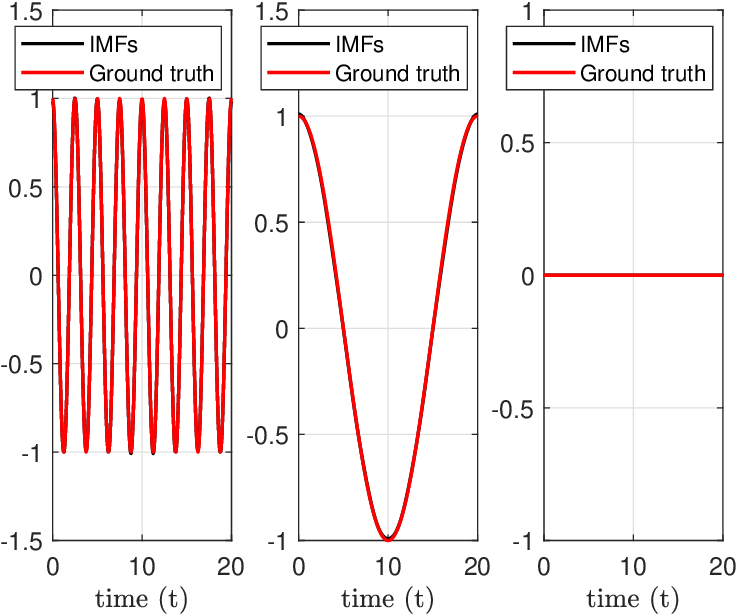}
    \caption{Example 1. Left panel: time decomposition for data at position $\mathbf{v}=(100,\,10)$. Right panel: time decomposition for data at position $\mathbf{v}=(150,\,150)$.}
    \label{fig:Ex1_IMFs_t}
\end{figure}

\begin{figure}[ht]
    \centering
\includegraphics[width=0.9\linewidth]{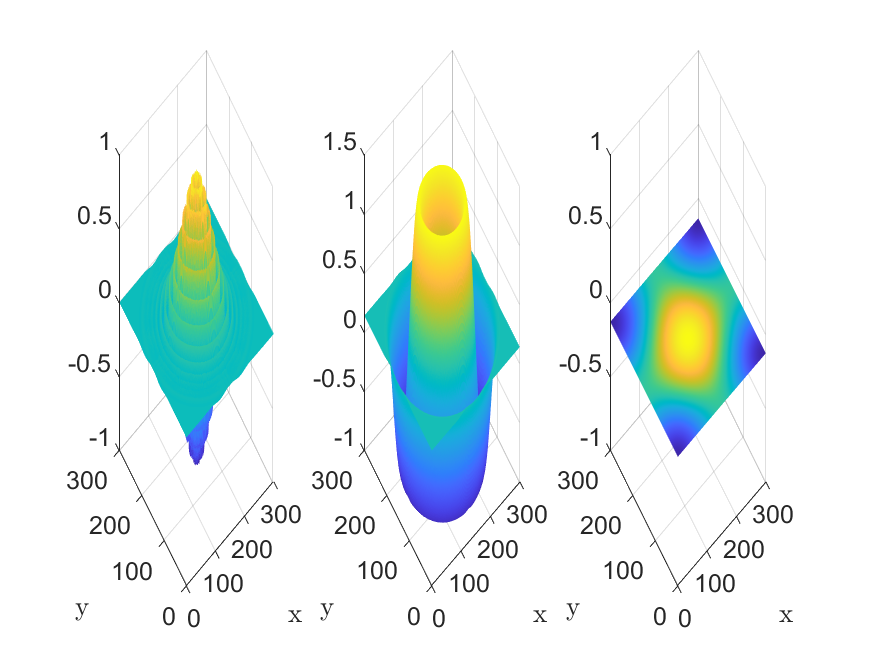}
    \caption{Example 1. From left to right, first, second IMF and trend over space.}
    \label{fig:Ex1_IMFs_s}
\end{figure}

\begin{figure}[ht]
    \centering
\includegraphics[width=0.9\linewidth]{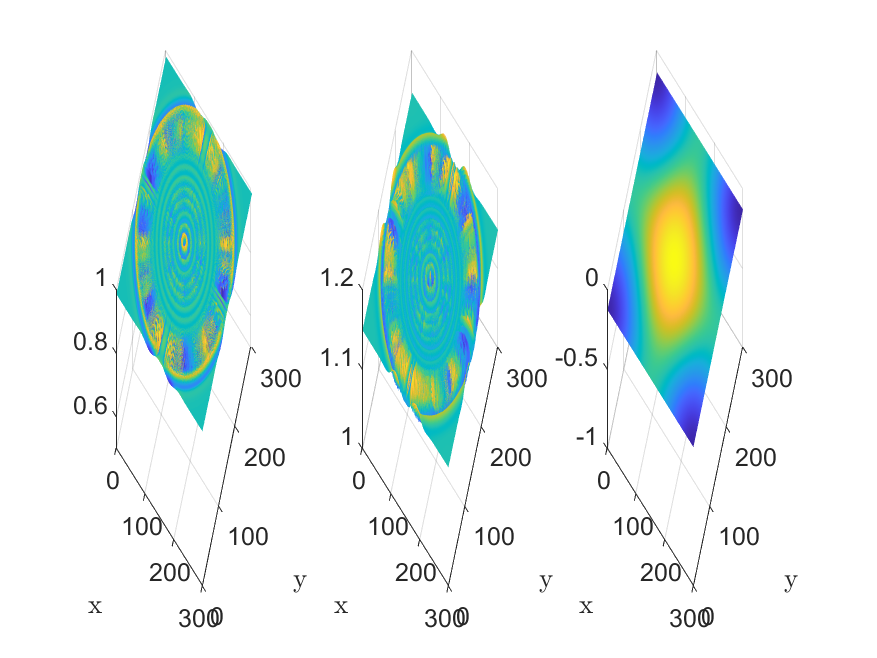}
    \caption{Example 1. Differences between the ground truth and the first, and second IMF and trend produced over space by MdMvFIF algorithm, respectively.}
    \label{fig:Ex1_IMFs_s_err}
\end{figure}

We apply MdMvFIF to decompose such signals into space and time IMFs. We plot the results in Figures \ref{fig:Ex1_IMFs_s} and \ref{fig:Ex1_IMFs_t}. In Figure \ref{fig:Ex1_IMFs_s_err} we report the differences between the ground truth space components and the IMFs produced with the proposed method. 

\subsection{Second artificial example}

In this second example, we have a signal containing two non-stationary frequencies active in space and two non-stationary frequencies active in time. 
\begin{figure}[ht]
    \centering
\includegraphics[width=0.9\linewidth]{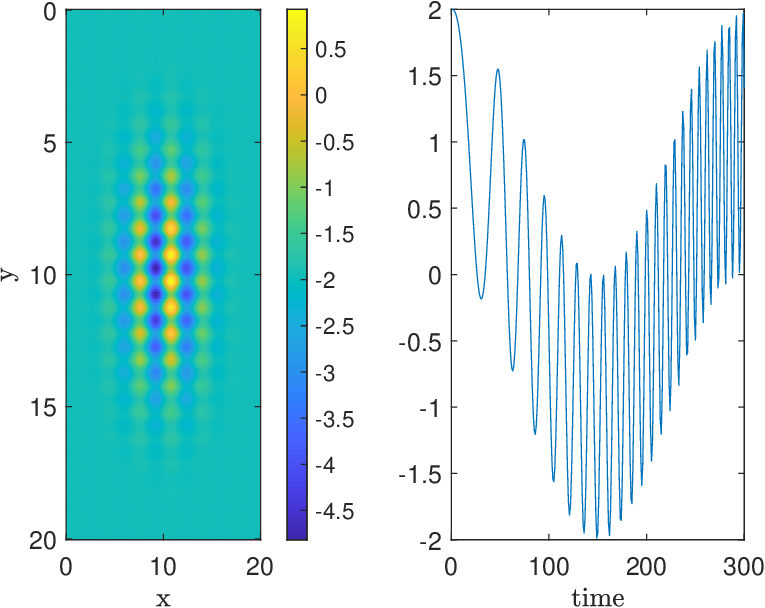}
    \caption{Example 2. Left panel: Signal at time $t=0$. Right panel: Time evolution for position $\mathbf{v}=(150,\,150)$.}
    \label{fig:Ex2_sig}
\end{figure}

\begin{figure}[ht]
    \centering
\includegraphics[width=0.9\linewidth]{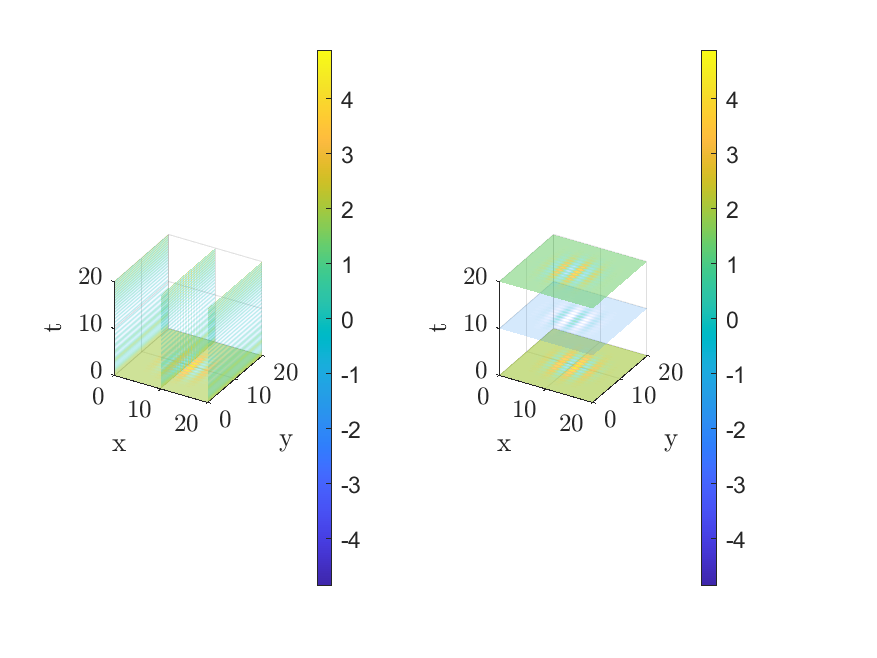}
    \caption{Example 2. Left panel: 3D signal evolution in the second space variable and time. Right panel: 3D signal evolution in space at different time stamps.}
    \label{fig:Ex2_sig2}
\end{figure}

The data are shown in Figures \ref{fig:Ex2_sig} and \ref{fig:Ex2_sig2}. 

\begin{figure}[ht]
    \centering
\includegraphics[width=0.45\linewidth]{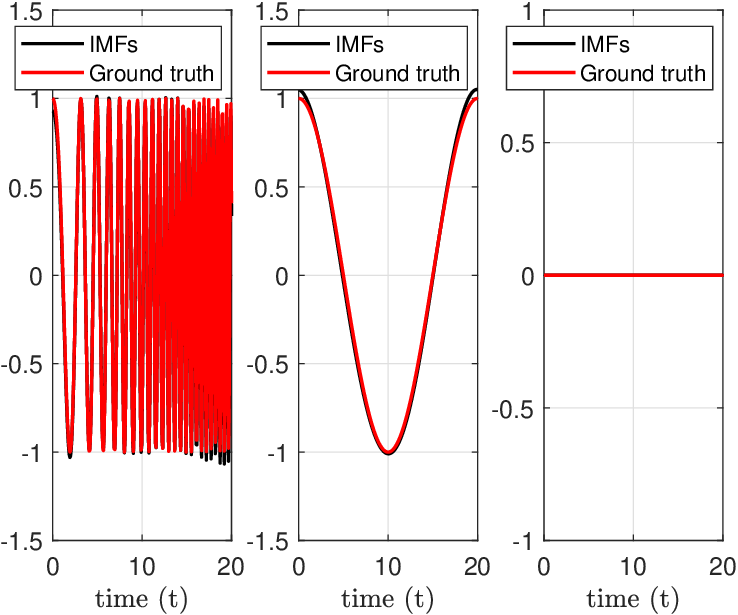}
\includegraphics[width=0.45\linewidth]{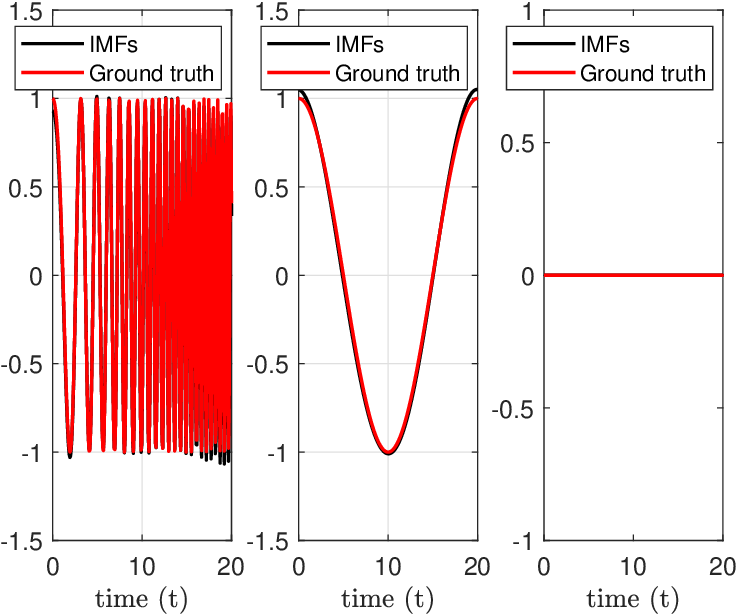}
    \caption{Example 2. Left panel: time decomposition for data at position $\mathbf{v}=(1,\,1)$. Right panel: time decomposition for data at position $\mathbf{v}=(100,\,10)$.}
    \label{fig:Ex2_IMFs_t}
\end{figure}

\begin{figure}[ht]
    \centering
\includegraphics[width=0.9\linewidth]{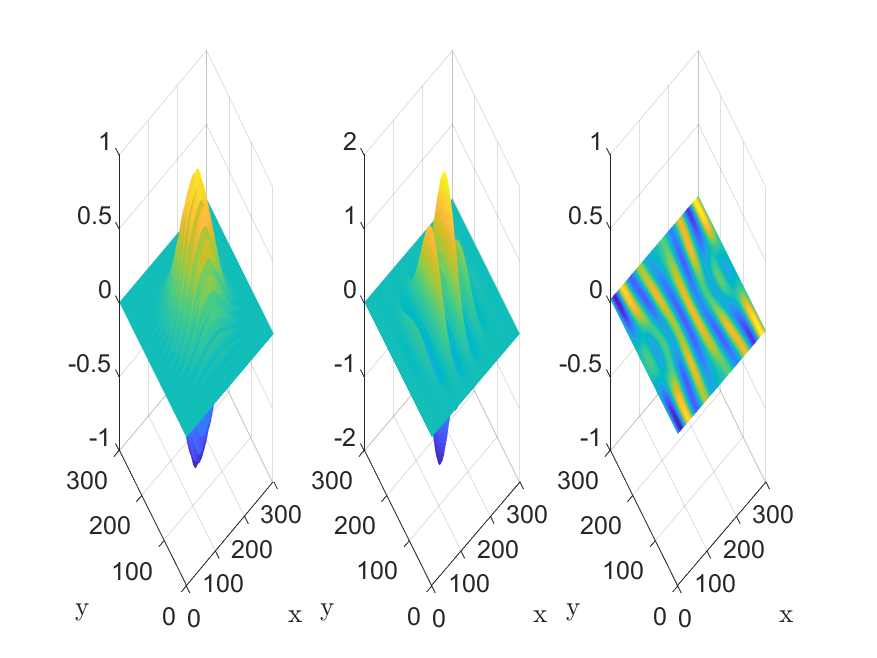}
    \caption{Example 2. From left to right, first, second IMF and trend over space.}
    \label{fig:Ex2_IMFs_s}
\end{figure}

\begin{figure}[ht]
    \centering
\includegraphics[width=0.9\linewidth]{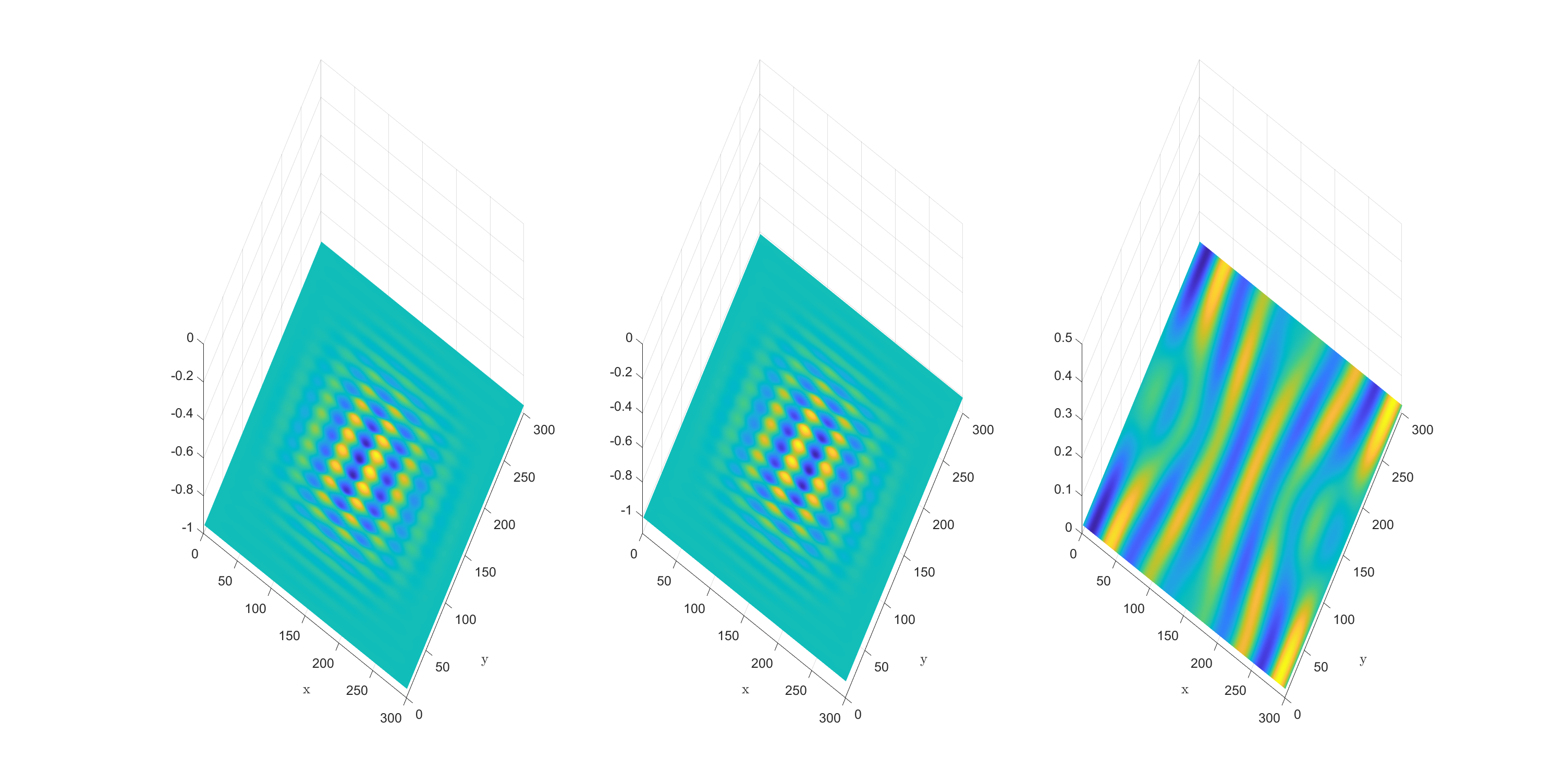}
    \caption{Example 2. Differences between the ground truth and the first, and second IMF and trend produced over space by MdMvFIF algorithm, respectively.}
    \label{fig:Ex2_IMFs_s_err}
\end{figure}

We apply MdMvFIF to decompose such signals into space and time IMFs. We plot the results in Figures \ref{fig:Ex2_IMFs_s} and \ref{fig:Ex2_IMFs_t}. In Figure \ref{fig:Ex2_IMFs_s_err} we report the differences between the ground truth space components and the IMFs produced with the proposed method. 

%


\subsection{Earth's air temperature}

We consider now the Earth's air temperature measured in kelvin on a global grid at an altitude of 2 meters from the surface of our planet. Such data set is made available through the NCEP/NCAR Reanalysis project\footnote{The original NCEP Reanalysis data have been provided by the NOAA/OAR/ESRL PSD, Boulder, Colorado, USA \url{http://www.esrl.noaa.gov/psd/} and have been edited by the Climatic Research Unit, University of East Anglia \url{https://crudata.uea.ac.uk/cru/data/ncep/}.} For more details on the project we refer the interested reader to \cite{kalnay1996ncep}.

\begin{figure}[ht]
    \centering
\includegraphics[width=0.9\linewidth]{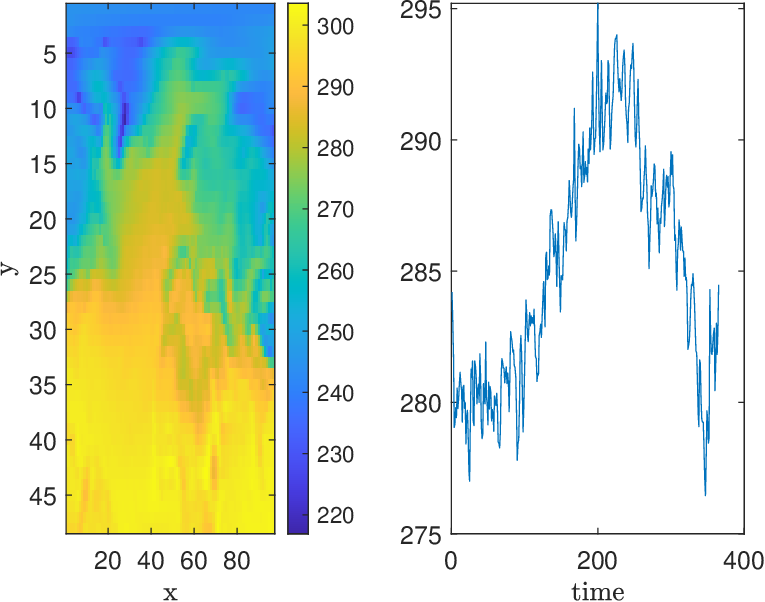}
    \caption{Example 3. Left panel: air temperature measured on January 1, 2022. Right panel: Time evolution for position $\mathbf{v}=(20,\,50)$.}
    \label{fig:Earth_sig}
\end{figure}

In this work, we consider the temperatures measured daily from January 1, 2022, to December 31\footnote{Data are available at \url{https://crudata.uea.ac.uk/cru/data/ncep/qs_eurasia/daily/sflux/air.2m/}}. A sample of the data is displayed in Figure \ref{fig:Earth_sig}.

The data $f(v_1,v_2,t)$ are stored as 3D tensors of dimension $48\times 97\times 365$. 

\begin{figure}[ht]
    \centering
\includegraphics[width=0.9\linewidth]{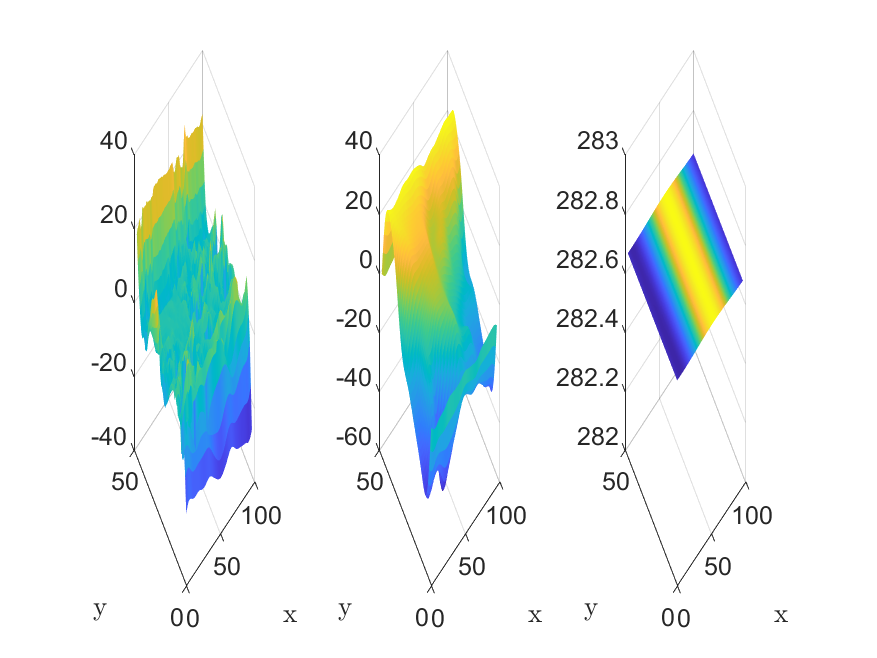}
    \caption{Example 3. Decomposition over space of the data into two IMFs plus a trend.}
    \label{fig:Earth_IMFs_s}
\end{figure}

\begin{figure}[ht]
    \centering
\includegraphics[width=0.45\linewidth]{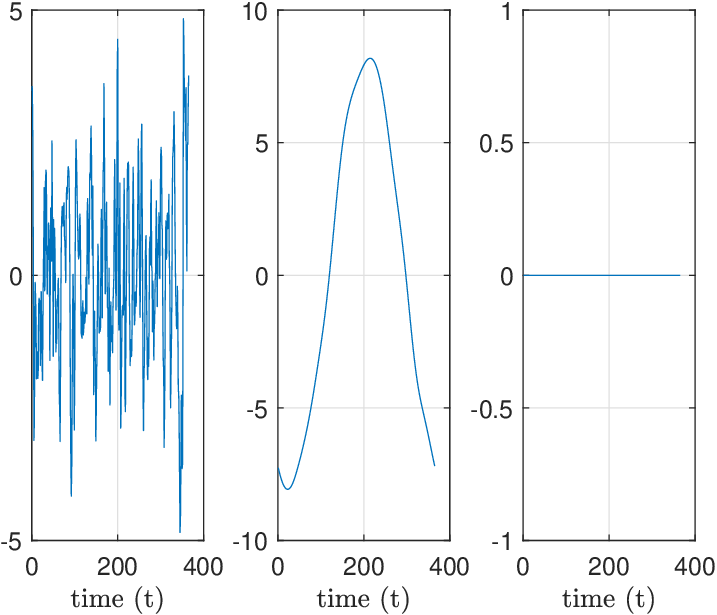}
\includegraphics[width=0.45\linewidth]{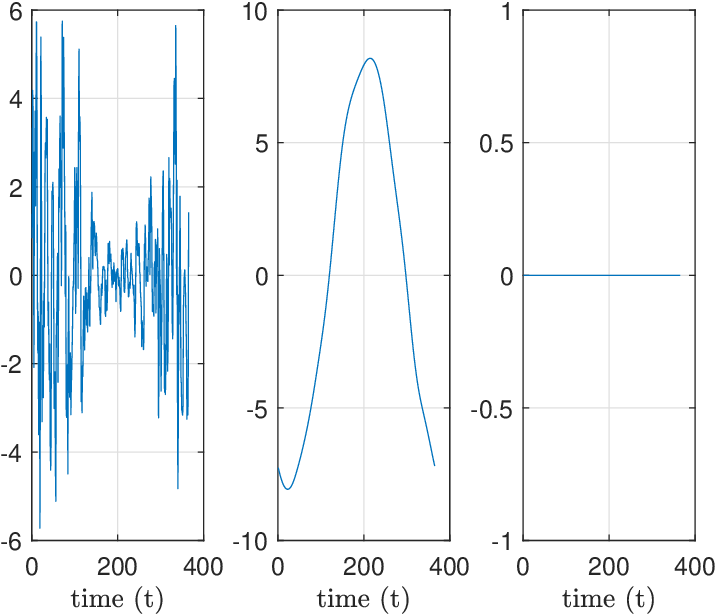}
    \caption{Example 3. Left panel: time decomposition for data at position $\mathbf{v}=(20,\,50)$. Right panel: time decomposition for data at position $\mathbf{v}=(10,\,45)$.}
    \label{fig:Earth_IMFs_t}
\end{figure}

We apply MdMvFIF to decompose such signals into space and time IMFs. We plot the results in Figures \ref{fig:Earth_IMFs_s} and \ref{fig:Earth_IMFs_t}. 

\section{Conclusion}\label{sec:Conclusion}

In this work, we introduced the Multidimensional and Multivariate Fast Iterative Filtering (MdMvFIF) method, an extension of the Fast Iterative Filtering (FIF) approach designed to decompose complex signals that vary both in space and time. This new algorithm overcomes several limitations present in previous signal decomposition methods, such as the separation of time and space variations and the inability to independently filter different frequency ranges in space and time.

We demonstrated that the MdMvFIF method can be applied to high-dimensional and multivariate datasets, handling not only multidimensional signals but also signals that evolve over time. By iterating the filtering process first in space and then in time, the MdMvFIF method effectively isolates Intrinsic Mode Functions (IMFs) from complex data, making it a valid tool for the study of many real-world data sets. We presented several numerical examples, including synthetic signals and geophysical data, to show the versatility and robustness of the approach. 

In particular, the MdMvFIF algorithm provides a powerful tool for signal analysis in applications where the interplay between spatial and temporal variations is complex and intertwined. Future work will explore further optimization of the algorithm, including strategies for speeding up calculations and its integration with machine learning techniques for feature extraction and analysis.

In conclusion, the MdMvFIF technique expands the capabilities of signal decomposition, offering a nonlinear, adaptive, and efficient approach for analyzing multidimensional and multivariate signals that vary over both space and time. This makes it a promising method for a wide range of applications, from geophysical data analysis to more general signal processing tasks.

\section*{Acknowledgment}
AC, EP and RC are members of the Gruppo Nazionale Calcolo Scientifico-Istituto Nazionale di Alta Matematica (GNCS-INdAM).

The research of AC and RC were partially supported through the GNCS-INdAM Project, CUP E53C22001930001 and CUP E53C23001670001. AC and EP were supported by the Italian Ministry of the University and Research and the European Union through the ``Next Generation EU'', Mission 4, Component 1, under the PRIN PNRR 2022 grant number CUP E53D23018040001 ERC field PE1 project P2022XME5P titled ``Circular Economy from the Mathematics for Signal Processing prospective'', and they thanks the Italian Ministry of the University and Research (MUR) for the financial support, CUP E13C24000350001, under the ``Grande Rilevanza'' Italy – China Science and Technology Cooperation Joint Project titled ``sCHans – Solar loading infrared thermography and deep learning teCHniques for the noninvAsive iNSpection of precious artifacts'', PGR02016.

\end{document}